\documentclass[a4paper]{article}
\usepackage{latexsym,a4,bbm}
\newtheorem{lem}{\bf Lemma}[section]
\newtheorem{de}[lem]{\bf Definition}
\newtheorem{thm}[lem]{Theorem}

\newenvironment{proof}{\noindent\emph{Proof.}}{\hspace*{\fill}{$\Box$}\bigskip\par}

\input{amssym.def}
\input{amssym}
\def\m{{\hfill\hbox to 0pt{$\rule{2pt}{3mm}$\hss}\par\noindent}}

\title{Construction of Diffusion Algebras}
\author{ P.N. Pyatov$^{a}$ and R. Twarock$^{b}$\\[\bigskipamount]
${}^{a}$ Bogoliubov Laboratory of Theoretical Physics, Joint Institute\\
for Nuclear Research, Dubna, Moscow Region 141980, Russia\\[\medskipamount]
${}^{b}$ Department of Mathematics, City University,\\
Northampton Square, London EC1V 0HB}
\begin{document}
\maketitle

\begin{abstract}
In \cite{IPR} Diffusion algebras have been introduced in the context of
one-dimensional stochastic processes with exclusion in statistical
mechanics. While this reference is focused on the needs of the physicist reader and thus states results without proofs and focuses on the discussion of lower-dimensional examples, it is the purpose of this paper to present a construction formalism for Diffusion algebras and to use the latter to prove the results in that reference.
\end{abstract}

\section{Introduction}

Diffusion algebras play a key role in the understanding of
one-dimensional stochastic processes. In the case of $N$ species of
particles with only nearest-neighbour interactions with exclusion on a
one-dimensional lattice, Diffusion algebras are useful tools in finding
expressions for the probability distribution of the stationary state of
these processes. Following the idea of matrix product states
\cite{MPS1,MPS2}, the latter are given in terms of monomials built from the
generators of a quadratic algebra. Depending on whether the system is
closed, i.e. the stochastic process is defined on a ring, or open, in which
case boundary conditions at the end of the lattice come into play, this
expression varies; \cite{IPR} presents an exposition of these facts and
the reader is referred to this reference and references within for more
details about the application of Diffusion algebras in physics.

It is the purpose of this work to treat Diffusion algebras from the
mathematician's point of view and to prove a construction theorem for Diffusion algebras.

\noindent We consider the following setting:

\noindent Let $\alpha < \beta \in I_N:=\lbrace 1, \ldots, N \rbrace$ and
consider quadratic relations of the form
\begin{equation}\label{algs}
g_{\alpha \beta}D_{\alpha}D_{\beta} - g_{\beta \alpha} D_{\beta} D_{\alpha} =
x_{\beta} D_{\alpha} - x_{\alpha} D_{\beta}
\end{equation}
with $g_{\alpha \beta}\in {\Bbb R}\setminus{\{0\}}$, $g_{\beta \alpha}\in {\Bbb R}$ and $x_{\alpha}\in {\Bbb
C}$~\footnote{
Note that dependence on all the nonvanishing coefficients $x_\alpha$
in (\ref{algs}) is easily suppressed by
rescaling the elements $D_\alpha$ as $x_\alpha
D_\alpha$.
We choose to display  the dependence on these
coefficients here,
because they are important in applications to stochastic models,
which is the physical motivation for the study of this type of algebras.
}.
Then on has
\begin{de}\label{Diff}
An algebra with generators $\lbrace D_{\alpha}\vert \alpha \in I_N\rbrace$
and relations of type
(\ref{algs}) is called \textbf{Diffusion algebra}, if it admits
a linear PBW-basis of ordered
monomials of the form
\begin{equation}\label{basis}
D_{\alpha_1}^{k_1}D_{\alpha_2}^{k_2}\ldots D_{\alpha_n}^{k_n}, \qquad k_j \in {\Bbb N}^0
\end{equation}
with ${\alpha_1}>{\alpha_2}>\ldots >{\alpha_n}$.
\end{de}

We remark that although we formulate the mathematical setting for  coefficients in ${\Bbb R}$ in the relations (\ref{algs}),
for physical reasons only relations with positive coefficients $g_{\alpha \beta}\in {\Bbb R}^{> 0}$ and  $g_{\beta \alpha}\in {\Bbb R}^{\geq 0}$ $(\alpha<\beta)$ are relevant, because they are interpreted as hopping rates in stochastic models.
Since we are treating Diffusion algebras from the mathematical point of view here, we will not impose this restriction, but comment on the implications of this restriction on our results after the main theorem.

The requirement of having a PBW basis (\ref{basis}) implies conditions on
the coefficients $g_{\alpha \beta}$ and $x_{\alpha}$ in (\ref{algs})
according to the Diamond Lemma in Ring Theory \cite{DL}. In particular, the
latter gives a criterion to check under which conditions the relations in
(\ref{algs}) are of PBW type: it is the case if each subset of three
generators $\lbrace D_\alpha, D_\beta, D_\gamma \rbrace$ with ordering
$\alpha < \beta < \gamma$ is reduction unique with respect to the ordering,
that is if the two ways of reducing the monomial $D_\alpha D_\beta D_\gamma$
to the monomial $D_\gamma  D_\beta D_\alpha$ lead to the same result when
expressed in the PBW basis (\ref{basis}).

The task of deriving all Diffusion algebras with $N$ generators thus reduces
to the following two steps:
\begin{enumerate}
\item Find all Diffusion algebras with three generators.
\item Find all algebras with $N$ generators such that each
subset of three generators coincides with one of the cases listed before.
\end{enumerate}

The first is a trivial exercise, which amounts to finding those coefficients
$g_{\alpha \beta}$ and $x_\alpha$ in (\ref{algs}) for which a set
$\{D_\gamma,  D_\beta, D_\alpha\}$ of three generators is reduction unique
in the above sense. The corresponding list of algebras is given in
\cite{IPR}, and we review it here in order to set up notation and render
this paper self-contained. The second is a combinatorial problem, and
requires to  combine in a consistent way the three generator algebras listed
before to algebras with $N$ generators for general $N > 3$.

A construction method for Diffusion algebras, and thus a constructive method
to approach the second point, is the so-called blending procedure in
the above reference, which is an inductive procedure for the construction of Diffusion algebras.
It uses the three generator cases and augments them to larger
units by attaching further generators in accordance with the requirements of
the Diamond Lemma, then giving a prescription how these larger building
blocks may be glued (or in the terminology of this reference ``blended'')
together in order to obtain a general Diffusion algebra of $N$ generators.
The advantage of the inductive procedure is that it facilitates the construction of representations, which are crucial for applications in physics. The main purpose of this paper is to provide a different construction method, which is more suitable for mathematical purposes, in particular, to deliver a proof for the fact that the
set of algebras obtained via the blending procedure corresponds exactly to
the set of Diffusion algebras in Definition \ref{Diff}.

After recalling the three generator case in Section 2, we
present in Section 3 a derivation of Diffusion algebras from first
principles. Furthermore, we present in this section a compact formulation for the blending procedure, and obtain the exhaustiveness of the inductive approach in \cite{IPR} as a corollary to our main theorem.

\section{Review of Diffusion algebras with 3 generators}

As mentioned in the introduction, the three generator case provides the
building blocks for the derivation of Diffusion algebras according to the
Diamond Lemma and we therefore briefly recall the results of \cite{IPR} for
this case.

Consider a set $\lbrace D_{\alpha}, D_{\beta}, D_{\gamma} \rbrace$ of three
generators with an ordering induced by the ordering of the index set $\alpha
< \beta < \gamma$ and relations as in (\ref{algs}). Since
$g_{\alpha\beta}\not = 0$ for all $\alpha<\beta\in I_N$ by assumption, we
can cast the relations into the following form
\begin{equation}\label{algs3}
\begin{array}{rcl}
D_{\alpha}D_{\beta} & = & q_{\beta \alpha}
 D_{\beta} D_{\alpha} + {x}^{\alpha\beta}_{\beta} D_{\alpha} - {x}^{\alpha\beta}_{\alpha} D_{\beta} \\
D_{\alpha}D_{\gamma} & = & q_{\gamma \alpha}
 D_{\gamma} D_{\alpha} + {x}^{\alpha\gamma}_{\gamma} D_{\alpha} - {x}^{\alpha\gamma}_{\alpha} D_{\gamma} \\
D_{\beta}D_{\gamma} & = & q_{\gamma \beta}
 D_{\gamma} D_{\beta} + {x}^{\beta\gamma}_{\gamma} D_{\beta} - {x}^{\beta\gamma}_{\beta} D_{\gamma} \,,
\end{array}
\end{equation}
where $q_{ji}:= \frac{g_{ji}}{g_{ij}}$, ${x}^{ij}_k:=\frac{x_k}{g_{ij}}$ for
$k$, $i<j\in\{\alpha, \beta, \gamma\}$. Then, using (\ref{algs3}), any
monomial can be expressed in terms of the PBW basis (\ref{basis}). This is
well defined, if applying (\ref{algs3}) in different orders leads to the
same result, that is if the  reductions
\begin{equation}
D_\alpha D_\beta D_\gamma \rightarrow D_\beta D_\alpha  D_\gamma \rightarrow D_\beta D_\gamma D_\alpha
\rightarrow D_\gamma D_\beta D_\alpha
\end{equation}
and
\begin{equation}
D_\alpha D_\beta D_\gamma \rightarrow  D_\alpha  D_\gamma D_\beta \rightarrow  D_\gamma D_\alpha D_\beta
\rightarrow D_\gamma D_\beta D_\alpha
\end{equation}
using (\ref{algs3}) coincide when expressed in the PBW basis (\ref{basis})
\cite{DL}. This leads to restrictions on the coefficients $g_{\alpha\beta}$ and
$x_{\alpha}$ in (\ref{algs}). In particular,
 they are
constrained by a set of 6 equations (see (2.5) --
(2.10) in \cite{IPR}) and their solutions determine all Diffusion algebras
of three generators. The latter are listed here for future convenience and
in order to set up notations. Throughout this section, we
assume $\alpha < \beta < \gamma$ and $x_j \not =0$ for
$j\in\{\alpha,\beta,\gamma\}$.

\begin{enumerate}

\item \underline{The case of $A_I$:}
\begin{equation}\label{first}
\begin{array}{rcl}
g [ D_\alpha , D_\beta ] & =  & x_\beta  D_\alpha -  x_\alpha D_\beta\\
g [ D_\alpha , D_\gamma ] & =  & x_\gamma  D_\alpha -  x_\alpha D_\gamma\\
g [ D_\beta , D_\gamma ] & =  & x_\gamma  D_\beta -  x_\beta D_\gamma
\end{array}
\end{equation}
where $g\neq 0$.

\item \underline{The case of $A_{II}$:}
\begin{equation}
\begin{array}{rcl}
g_{\alpha\beta} D_\alpha D_\beta  & = & x_\beta  D_\alpha -  x_\alpha  D_\beta\\
g_{\alpha\gamma} D_\alpha D_\gamma  & = & x_\gamma  D_\alpha -  x_\alpha  D_\gamma\\
g_{\beta\gamma} D_\beta D_\gamma  & = & x_\gamma  D_\beta -  x_\beta  D_\gamma
\end{array}
\end{equation}
where  $g_{ij}:= g_i - g_j$ with $ g_i\neq g_j$ for all $i<j\in\{\alpha,\beta,\gamma\}$.

\item \underline{The case of $B^{(1)}$:}
\begin{equation}\label{B1}
\begin{array}{rcl}
g_{\beta} D_{\alpha} D_{\beta} - (g_{\beta}-\Lambda) D_{\beta} D_{\alpha} &=& -x_{\alpha}D_\beta  \\
g  D_{\alpha} D_{\gamma}- (g-\Lambda) D_{\gamma} D_{\alpha} &=& x_\gamma D_\alpha -  x_\alpha D_\gamma \\
g_{\beta} D_{\beta} D_{\gamma} -(g_{\beta}-\Lambda) D_{\gamma} D_{\beta} &=& x_\gamma  D_\beta
\end{array}
\end{equation}
where $g\neq 0$ and $g_\beta\neq 0$.
For the
same ordering, we also find relations of type $B^{(1)}$ which are relations
(\ref{B1}) with an exchange $\alpha \leftrightarrow \beta$ or $\gamma
\leftrightarrow \beta$ and restrictions $g\neq 0$ \mbox{and}
$g_\alpha\notin \{0,\Lambda\}$
or, respectively, $g\neq 0$ and
$g_\gamma \notin \{0,\Lambda\}$ on the parameters.

\item \underline{The case of $B^{(2)}$:}
\begin{equation}
\begin{array}{rcl}
g_{\alpha \beta} D_{\alpha} D_{\beta} & =& -  x_\alpha D_\beta \\
g_{\alpha \gamma}  D_{\alpha} D_{\gamma} -g_{\gamma \alpha} D_{\gamma} D_{\alpha} & = &
x_\gamma D_\alpha -  x_\alpha D_\gamma \\
g_{\beta \gamma} D_{\beta} D_{\gamma}  & =& x_\gamma  D_\beta
\end{array}
\end{equation}
where $g_{\alpha\beta}$, $ g_{\alpha\gamma}$ and $g_{\beta\gamma}\neq 0$.

\item \underline{The case of $B^{(3)}$:}
\begin{equation}\label{B3}
\begin{array}{rcl}
g D_{\alpha} D_{\beta}  -(g-\Lambda) D_{\beta} D_{\alpha} & =&
x_\beta D_\alpha -  x_\alpha D_\beta \\
g_\gamma  D_{\alpha} D_{\gamma} & = &  -x_\alpha D_\gamma \\
(g_{\gamma}-\Lambda) D_{\beta} D_{\gamma} & = & -x_\beta D_\gamma
\end{array}
\end{equation}
where $g\neq 0$ and $g_\gamma\notin\{0,\Lambda\}$.

\item \underline{The case of $B^{(4)}$:}
\begin{equation}\label{B4}
\begin{array}{rcl}
(g_{\alpha}-\Lambda) D_{\alpha} D_{\beta} & =& x_\beta D_\alpha \\
g_{\alpha}  D_{\alpha} D_{\gamma} & =& x_\gamma D_\alpha  \\
g D_{\beta} D_{\gamma} - (g-\Lambda) D_{\gamma} D_{\beta} & =& x_\gamma D_\beta - x_\beta D_\gamma
\end{array}
\end{equation}
where $g\neq 0$ and $g_\alpha\notin\{0,\Lambda\}$.

\item \underline{The case of $C^{(1)}$:}
\begin{equation}\label{C1}
\begin{array}{rcl}
g_\beta D_{\alpha} D_{\beta} - (g_\beta-\Lambda) D_{\beta} D_{\alpha} & =& - x_\alpha D_\beta \\
g_\gamma D_{\alpha} D_{\gamma} - (g_\gamma-\Lambda) D_{\gamma} D_{\alpha} & = &
 - x_\alpha D_\gamma \\
g_{\beta \gamma} D_{\beta} D_{\gamma} - g_{\gamma \beta} D_{\gamma} D_{\beta} & = & 0
\end{array}
\end{equation}
where $g_\beta, g_\gamma$ and $g_{\beta\gamma} \neq 0$.
For the same ordering, we also
find relations of type $C^{(1)}$ which are relations (\ref{C1}) with an
exchange $\alpha \leftrightarrow \beta$
or $\alpha \rightarrow \gamma
\rightarrow \beta \rightarrow \alpha$
and restrictions $g_\alpha\neq\Lambda$ and $g_\gamma , g_{\alpha\gamma} \neq 0$
or, respectively,
$g_\alpha, g_\beta\neq\Lambda$ and
$g_{\alpha\beta} \neq 0$ on the parameters.

\item \underline{The case of $C^{(2)}$:}
\begin{equation}\label{C2}
\begin{array}{rcl}
g_{\alpha \beta} D_{\alpha} D_{\beta} -
g_{\beta \alpha} D_{\beta} D_{\alpha} & =& -  x_\alpha D_\beta \\
g_{\alpha \gamma} D_{\alpha} D_{\gamma} -
g_{\gamma \alpha} D_{\gamma} D_{\alpha} & =& -  x_\alpha D_\gamma \\
D_{\beta} D_{\gamma} & =& 0
\end{array}
\end{equation}
where $g_{\alpha\beta}$ and $g_{\alpha\gamma}\neq 0$. For the
same ordering, we also find relations of type $C^{(2)}$ which are relations
(\ref{C2}) with an exchange $\alpha \leftrightarrow \beta$
or $\alpha
\rightarrow \gamma \rightarrow \beta \rightarrow \alpha$
and
restrictions $g_{\alpha\beta}, g_{\beta\gamma}\neq 0$ or, respectively,
$g_{\alpha\gamma}, g_{\beta\gamma}\neq 0$ on the parameters.

\item \underline{The case of $D$:}
With $q_{ji}:=\frac{g_{ji}}{g_{ij}}$, $i$, $j\in\{\alpha,\beta,\gamma\}$
(recall that $g_{ij}\not = 0$ for $i < j$)  we have
\begin{equation}\label{last}
\begin{array}{rcl}
D_\alpha D_\beta - q_{\beta\alpha} D_\beta D_\alpha & = & 0\\
D_\alpha D_\gamma - q_{\gamma\alpha} D_\gamma D_\alpha & = & 0\\
D_\beta D_\gamma - q_{\gamma\beta} D_\gamma D_\beta & = & 0
\end{array}\,.
\end{equation}
\end{enumerate}

We remark that the division into algebras of type $A$, $B$, $C$ and $D$
reflects the number of coefficients $x_j$, $j\in \{ \alpha, \beta, \gamma
\}$, being zero in comparison with the general form (\ref{algs}): for
algebras of type $A$, $B$, $C$ and $D$ none respectively one, two, or all
three of the coefficients $x_i$ vanish. The subdivision for each type then
corresponds to the different choices for the coefficients $g_{\alpha \beta}$
which are compatible with the Diamond Lemma.

\section{The case of general $N$}

This section consists of four parts: we start by providing a decomposition
of the index set which later facilitates the presentation of the algebras. In other words, we decompose the whole family of algebras, which depends on the
ordered set of parameters $\{g_{\alpha\beta},x_\alpha\vert\alpha,\beta\in I\}$ in the relations (\ref{algs}), into several subfamilies. Each subfamily is determined by a specific subset of the parameters $x_\alpha$ and $g_{\alpha\beta}$, which are subject to a set of conditions formulated below (see conditions (\ref{Co1}), (\ref{requ1})--(\ref{extra2}), (\ref{Co2})
and (\ref{splitT}) below).

As a next step, we list some general properties  specific to Diffusion
algebras in each of the subfamilies. They are later used in the proof of the
main result. This is followed by the list of Diffusion algebras
and a theorem which proves the exhaustiveness of the approach.
We finally comment on the counting of Diffusion algebras.

\subsection{Decomposition of the index set}

The structure of the algebras in (\ref{first}) -- (\ref{last}) suggests the
following decomposition of the index set $I_N=\{1,\ldots,N\}$:
\begin{equation}
I_N = I \cup R
\end{equation}
where
\begin{equation}
\label{Co1}
\begin{array}{rcl}
I & := & \{ \alpha\in I_N | x_{\alpha} \neq 0 \}\\
R & := & \{ \alpha\in I_N | x_{\alpha} = 0 \}\,.
\end{array}
\end{equation}
We will use in the following the notation $N_I := \vert I \vert$ and $N_R :=
\vert R \vert$ for the cardinalities of these sets.

We introduce the following terminology and notations:

\begin{de}\label{normalOrdering}
Normal ordering of two generators $D_\alpha$ and $D_\beta$ is defined as
\begin{equation}\label{normal}
:D_\alpha D_\beta:\quad  := \left
\lbrace
\begin{array}{rcl}
D_\alpha D_\beta &  \mbox{ if } & \alpha<\beta \\
D_\beta D_\alpha &  \mbox{ if } & \beta<\alpha
\end{array}
\right.
\end{equation}
\end{de}
\begin{de}
For $\alpha < \beta$ we introduce the following  short-hand notation:
\begin{equation}
[D_{\alpha},D_{\beta}]_{q_{\beta\alpha}}:=D_{\alpha} D_{\beta} - q_{\beta\alpha} D_{\beta} D_{\alpha}
\end{equation}
where the index at the commutator is referring to the coefficients $q_{\beta
\alpha}$ in terms of which the commutator is defined.
\end{de}

Using these notation, we subdivide the set $R$
into nonitersecting and
non-empty subsets
\begin{equation}\label{decom}
R  :=  R_1 \cup R_2 \cup \ldots \cup R_{M_R}
\end{equation}
according to the following requirements:
\begin{itemize}
\item Relations between generators from the sets $R_a$ and $R_b$ for $a\not= b$ are given by:
\begin{equation}\label{requ1}
 :D_{r_1}D_{r_2}:\ = 0 \qquad \forall {r_1} \in R_a \mbox{~and~} \forall{r_2} \in
R_b.
\end{equation}
\item Relations within a set $R_a$  such that
 $\vert R_a\vert \geq 2$ are given by:
\begin{equation}\label{requ2}
{}[D_{r_1}, D_{r_2}]_{q_{r_2r_1}} =0 \qquad  \forall {r_1}<{r_2} \in R_a\ ,
\end{equation}
where the coefficients in (\ref{requ2}) are subject
to the condition opposite to (\ref{requ1}),
that is: for any subdivision
$R_a= R'\cup R''$ into two nonintersecting and
non-empty parts $R'$ and $R''$
\begin{equation}\label{extra}
\exists
{r_1}\in R'  \mbox{~and~} {r_2}\in R''\ :\
g_{r_1 r_2}g_{r_2r_1}\not = 0\,.
\end{equation}
In other words, this means that for any pair of indices $r, s\in R_a$ there
exists a finite sequence $\{r_k\in R_a\vert k=1,\dots ,n\}$ such that
$r_1=r$, $r_n=s$ and
\begin{equation}\label{extra2}
\prod_{k=1}^{n-1} g_{r_k r_{k+1}}g_{r_{k+1}r_k}\not = 0\ .
\end{equation}
Thus, the relations (\ref{extra}) and (\ref{extra2}) may be represented graphically via a {\em connectivity condition} on an ordered graph the vertices of which are labelled by the indices
$r\in R_a$ and the edges connect only those vertices $r_1<r_2$ for which the
condition $q_{r_2r_1}\neq 0$ is satisfied.
\end{itemize}

\noindent
Furthermore, for $N_I\geq 2$ we split the set $R$ into two sets $S$ and $T$ as follows:
\smallskip

For any $R_a\subset R$ we define
\begin{equation}\label{Co2}
R_a:=\left\{
\begin{array}{rl}
S_a & \mbox{ if }\exists r\in R_a\ \mbox{and} \, i\in I: \, g_{ir}g_{ri}\not = 0\ ,
\\[1mm]
T_a & \mbox{ otherwise }\,.
\end{array}
\right.
\end{equation}

Suppose that the $M_R$ sets $R_a$ in (\ref{decom}) split into $M_S$ sets $S_a$ and $M_T$ sets
$T_a$ in this way,  thus $M_R=M_S+M_T$. We number these sets as
$S_a$, $a=1,\ldots,M_S$,
and $T_a$, $a=1,\ldots,M_T$, and introduce
\begin{equation}\label{ST}
S \ := \ \cup_{a=1}^{M_S} S_a\ , \qquad
T \ := \ \cup_{a=1}^{M_T} T_a  \,.
\end{equation}
Although the decomposition of the set $S$ into subsets $S_a$ has been used in the definition of the set $S$, it will not be of practical importance
in what follows. Contrary to that,
the structure of the set $T$ is crucial and
needs further refinement.
\smallskip

For any $T_a\subset T$  define
\begin{equation}\label{splitT}
T_a:=\left\lbrace
\begin{array}{rlc}
T_a^{\bullet} & \mbox{ if~}
\exists\ i<j\in I\ :& T_a\subset\{i+1,i+2\dots ,j-1\}
\mbox{~and~}
\\[1mm] & &
I\cap \{i+1,i+2,\dots ,j-1\} =\emptyset\ ,
\\[1mm]
T_a^{\circ} & \mbox{ otherwise }\,.&
\end{array}
\right.
\end{equation}
Thus in short hand notation $T= \{T_a^{\bullet}\vert a=1,\ldots,M_T^{\bullet}\}\cup
\{T_a^{\circ}\vert a=1,\ldots,M_T^{\circ}\}$
with $M_T=M_T^{\bullet}+M_T^{\circ}$.

\subsection{General structural remarks about $N$-generator Diffusion algebras}

Until now we have primarily discussed index sets.
By an abuse of terminology, we will from now on also refer to ``generators of a set $I$, $S$, $T$, or $R$'' meaning the generators indexed by elements from the corresponding set.

\begin{de}
A set of three generators $\{ D_x, D_y, D_z \}$ with $x$, $y$ and $z\in X$, $Y$, $Z$, respectively, where
$X$, $Y$ and $Z$ are any of the sets $I$, $R$, $S$ and $T$ or any set in their decomposition will be
called a triplet (of type) $\{X,Y,Z\}$.
\end{de}

Note that any triplet of type $\{I,I,I\}$ in a Diffusion algebra of $N\geq 3$
generators gives rise to a Diffusion algebra of type $A_I$ or $A_{II}$, any triplet of
type $\{I,I,R\}$ to a Diffusion
algebra of type $B^{(1)}$, $B^{(2)}$, $B^{(3)}$, or $B^{(4)}$, any triplet of type
$\{I,R,R\}$ to a Diffusion algebra of type $C^{(1)}$ or $C^{(2)}$ and any triplet of type $\{R,R,R\}$
to a Diffusion algebra of type $D$.

Then we have:

\begin{lem}\label{lem}\smallskip
For any Diffusion algebra (\ref{algs}) with $N\geq 3$ generators the following
statements hold
\begin{enumerate}
\item
If $N_I \geq 3$, then all subalgebras corresponding to triplets of type $\{I,I,I\}$ are
of the same type, which is either  $A_I$ (that is, $g_{ij}=g$ $\forall i,j \in I$) or $A_{II}$
(that is, $g_{ji}=0$, $g_{ij}=g_i-g_j$,
$g_i\neq g_j$, $\forall i<j \in I$).
\item
If $N_I \geq 3$ and
all subalgebras corresponding to triplets $\{I,I,I\}$ are of
type $A_I$ then for any
$s\in S$ and for all $i\in I$ one has
\begin{equation}
\label{vot}
g_{is}=g_{si}=g_s\ .
\end{equation}
\item
If $N_I \geq 3$ and
all subalgebras corresponding to triplets $\{I,I,I\}$ are of
type $A_{II}$  then $S=\emptyset$.
\item
Let $N_I\geq 2$.

For any $i\in I$ and for all $t, t'\in T_a$ (here $T_a$ means both
$T^\circ_a$ and $T^\bullet_a$) with $t<i$ and $t'>i$
the
coefficients $g_{t i}$ and $g_{it'}$
depend only on the index $a$ of the set $T_a$ and not on the individual
indices $t$ or $t'$. If $t, t'\in T^\circ_a$ one furthermore has
$g_{t i}=-g_{i t'}$.

For any $i<j\in I$ and for all $t,t'\in T_a$: $t<i$ and $t'>j$
\begin{eqnarray}
\label{tij}
g_{t i} + \Lambda_{ij} &=& g_{t j}\ ,
\\
\label{ijt}
g_{i t'} &=& g_{j t'} +\Lambda_{ij}\ ,
\end{eqnarray}
where $\Lambda_{ij}:= g_{ij}-g_{ji}$.
\item
Let $N_I = 1$. Denote the only index in $I$ as $\mathbf i$ in order to
stress that it is not a running index. For all $r\in R_a$
one has
\begin{equation}
\label{ri}
g_{{\mathbf i}r}- g_{r {\mathbf i}} = \Lambda_a\ .
\end{equation}
Note that both the left and the right hand sides of relation (\ref{ri}) depend only on the
index $a$ of the set $R_a$ and not on the individual index $r$.
\end{enumerate}
\end{lem}
\begin{proof}
{\bf 1.~}
It follows from the fact that in each set $\{D_i, D_j, D_k| i,j,k \in I\}$
one has either $g_{\beta \alpha}\neq 0$ for all $\alpha< \beta\in
\{i,j,k\}$, or, $g_{\beta \alpha}= 0$, for all $\alpha< \beta\in
\{i,j,k\}$, but no mixture thereof, which contradicts a mixing of $A_I$ and
$A_{II}$ type algebras.
\medskip

\noindent
{\bf 2.~}
Let $N_I \geq 3$ and $g_{ij}=g$ $\forall i,j \in I$. By (\ref{ST}) it is
enough to check (\ref{vot}) for
any $S_a\subset S$. Consider an index $r\in S_a$
which satisfies the condition $g_{rj}g_{jr}\neq 0$
for some $j\in I$. Then,
for any $i\in I$, the triplets $\{D_r,D_i,D_j\}$ are all of type
$B^{(1)}$ with $\Lambda=0$ and, hence, $g_{ir}=g_{ri}=g_r$.

Next, take any $s\in S_a$.
By definition, there exists a sequence $r_k\in S_a$,
$k=1,\dots ,n$, such that
$r_1=r$, $r_n=s$ and such that the connectivity condition (\ref{extra2})
is satisfied.
Then, for any $i\in I$, starting with the triplet $\{D_{r_1},D_{r_2},D_i\}$
one inductively proves that all the triplets $\{D_{r_k},D_{r_{k+1}},D_i\}$
are of type
$C^{(1)}$ with $\Lambda=0$, and hence (\ref{vot}) follows.
\medskip

\noindent
{\bf 3.~}
Let $N_I \geq 3$ and $g_{ji}=0$, $g_{ij}=g_i-g_j$, $g_i\neq g_j$, $\forall i<j \in I$ and suppose
$S\neq\emptyset$. Consider some
$S_a\subset S$ and take those indices $r\in S_a$ and $i_0\in I$
for which
condition $g_{i_0r}g_{ri_0}\neq 0$ is satisfied. For any $j,k\in I$ the triplets
$\{D_{i_0}, D_j, D_r\}$ and $\{D_{i_0}, D_k, D_r\}$ are both of  type $B^{(1)}$,
which also implies
that the triplets $\{D_j,D_k,D_r\}$, $\forall j,k\in I$, are all of  type $B^{(1)}$. Now, there is
no mutual ordering of any arbitrarily chosen indices $i<j<k\in I$
and the index $r\in S_p$ for which the existence of any $B^{(1)}$-type triplet
$\{D_i,D_j,D_r\}$, $\{D_i,D_k,D_r\}$ and $\{D_j,D_k,D_r\}$ is compatible with the
condition that $g_{ik}=g_{ij}+g_{jk}$ --- a contradiction.
\medskip

\noindent
{\bf 4.~}
Let $N_I\geq 2$ and consider any three indices $i\in I$ and $t<t'\in T_a$.
Exploiting the connectivity property (\ref{extra2}) of the set $T_a$
one can find a sequence $\{t_k\vert k=1,\dots ,n\}$
such that $t_1=t$, $t_n=t'$ and such that all
the $C$-type triplets $\{D_i,D_{t_k},D_{t_{k+1}}\}$ are not of  type $C^{(2)}$.
Hence, their corresponding nonzero coefficients $g_{it_k}$ (for $i<t_k$) or $g_{t_k i}$ (for $t_k<i$)
are subject to the conditions for triplets of type $C^{(1)}$ (see (\ref{C1}))
which
together with the definition (\ref{Co2}) of the set $T_a$
implies $g_{it}=g_{it'}$ in the case $i<t<t'$, $g_{ti}=g_{t'i}$ in the case $t<t'<i$
and $g_{ti}=-g_{it'}$ in the case $t<i<t'$, thus proving the first part of the
fourth statement in the Lemma.

To prove the second part notice that for any four indices $i,j\in I$ and $t,t'\in T_a$
which are ordered as $t<i<j<t'$ their corresponding triplets
$\{D_t,D_i,D_j\}$ and $\{D_i,D_j,D_{t'}\}$ are of type
$B^{(4)}$ and, respectively, $B^{(3)}$. Conditions (\ref{tij}) and (\ref{ijt})
then reproduce the relations between the coefficients in those triplets
(see eqs.(\ref{B3}) and (\ref{B4})).
\medskip

\noindent
{\bf 5.~}
Let $N_I=1$ and consider any pair of indices $r, r'\in R_a$.
As before, for every connective set $R_a$
there exists a chain of $C$-type triplets $\{D_{\mathbf i},D_{r_k},D_{r_{k+1}}\}$,
$k=1,\dots ,n-1$, with $r_1=r$ and $r_n=r'$ which are not of type $C^{(2)}$. Hence one obtains (\ref{ri}) with one and the same coefficient $\Lambda_a$ for all
$C^{(1)}$ type  triplets $\{I,R_a,R_a\}$.
\end{proof}

Lemma \ref{lem} suggests to list Diffusion algebras in families according
to the number of generators in the set $I$ and
provides information about the structure of relations among generators
from the sets $I$, $S$, $T$ and $R$ in each case.

\subsection{List of Diffusion algebras with $N$ generators}

In this subsection, we list all $N$-generator Diffusion algebras and provide
a theorem which proves the exhaustiveness of the formalism.

Diffusion algebras with $N$ generators are listed as five families of
algebras:
 $A_I$, $A_{II}$, $B$, $C$ and $D$.
As in the case of $N=3$
the number of nonzero coefficients $x_\alpha$ or, in other words, the
cardinality of the set $I$
is used as a criterion for separating Diffusion algebras into families of
the types $A$($N_I\geq 3$), $B$($N_I=2$), $C$($N_I=1$) or $D$($N_I=0$).
Type $A$ algebras are  separated further into two families
$A_{I}$ and $A_{II}$ depending on the number of nonzero coefficients
$g_{ij}$ with indices $i,j$ in the set $I$.

Different algebras in the families are obtained in dependence on the
choice
 of the decomposition of the set $I_N=\{1,2,\dots ,N\}$ into {\em ordered
subsets}
$I$, $S$, $T_a^{\circ}$, $a=1,\dots ,M_T^\circ$,
$T_b^{\bullet}$, $b=1,\dots ,M_T^\bullet$ (or $R_a$, $a=1,\dots M_R$ for $N_I=1$)
as well as on the choice of coefficients  in their defining relations.
Below we adopt a notation for Diffusion algebras where
the corresponding decomposition of the set $I_N$ is given explicitly
as argument to the family symbol.
To avoid any confusion let us stress that subscript
indices $a$ and $b$ in  our notation
 are treated as running ones so that, e.g.,
$$
A_I(I,S,T^\circ_a,T^\bullet_b)\equiv
A_I(I,S,T^\circ_1,\dots, T^\circ_{M^\circ_T},T^\bullet_1,\dots ,
T^\bullet_{M^\bullet_T})\ ,
$$
where we imply $I_N = I\cup S\cup (\cup_{a=1}^{M^\circ_T}T^\circ_a)\cup
(\cup_{b=1}^{M^\bullet_T}T^\bullet_b)$, and $I$, $S$, $T^\circ_a$ and
$T^\bullet_b$ are mutually nonintersecting ordered subsets in $I_N$.
The values
of the coefficients $g_{\alpha\beta}$ are not shown explicitly in these
notations so that in fact our notation displays  connnective components in a variety of Diffusion
algebras rather than the particular algebras.

All relations in (\ref{firstOne}) --
(\ref{lastOne}) below are to be complemented by relations (\ref{requ1}),
(\ref{requ2}) for the elements of the subset $R$
together with the conditions (\ref{extra}) or (\ref{extra2}) on the
coefficients involved.

\medskip
\noindent
1.
\underline{Diffusion algebras of type $A_I(I,S,T^\circ_a, T^\bullet_b)$,
 ~$N_I \geq 3$}:

\begin{eqnarray}
g[D_i,D_j] &  = & x_j D_i- x_i D_j \,,
\quad \forall i,j\in I, \nonumber\\
g_s [D_s, D_i] & = & x_i D_s \,,
\qquad \qquad \,\, \forall s \in S,\ i \in I, \nonumber\\
g_a^{\circ} :D_i D_{t}: & = &   - x_i D_{t}\,,
\qquad \;\quad \,\,  \forall a,\ t\in T_a^{\circ},\ i \in I,\label{firstOne}\\
g_{b}^{+} D_i D_{t} & = &  - x_i D_{t} \,,
\qquad \;\quad \,\,  \forall b,\ t\in T_b^{\bullet},\ i \in I : i<t, \nonumber\\
g_{b}^{-}  D_{t} D_i & = & x_i D_{t}\,,
\qquad \qquad \,\, \forall b,\ t\in T_b^{\bullet},\ i \in I,  :i>t, \nonumber
\end{eqnarray}
where $g$, $g_s$, $g_a^\circ$, $g_b^\pm\neq 0$.
\bigskip

\noindent
2.
\underline{Diffusion algebras of type $A_{II}(I,T^\circ_a, T^\bullet_b)$,
$N_I \geq 3$:}

\begin{eqnarray}\label{AII}
(g_{i}-g_j) D_iD_j & = & x_j D_i-x_i D_j \,,
\quad \forall i<j\in I,  \nonumber\\
(g_i+g_a^\circ) :D_i D_{t}: & = & -x_i D_{t} \,,
\qquad \quad \,\,\, \forall a,\ t\in T_a^\circ,\ i \in I, \\
(g_i+g_b^+) D_i D_{t} & = & -x_i D_{t} \,,   \qquad \quad \,\,\,
\forall b,\ t\in T_b^\bullet,\ i \in I : i< t,  \nonumber\\
(g_b^--g_i) D_t D_i  & = & x_i D_{t} \,,
\qquad \qquad \,\, \forall b,\ t\in T_b^\bullet,\ i \in I : i>t, \nonumber
\end{eqnarray}
where $g_i\neq g_j$ for $i\neq j$ and
$g_i\not\in\{g^\circ_a,\mp g^\pm_b\}$.
\bigskip

\noindent
3.
\underline{Diffusion algebras of type
$B(I=\{{\mathbf i},{\mathbf j}\},S,T^\circ_a, T^\bullet_b)$:}
\bigskip

We use the notation $\mathbf i$ and $\mathbf j$ with
${\mathbf i}<{\mathbf j}$ for the two elements of
the set $I$ to emphasize that they are not running indices.
Note also that ${\mathbf i}<t<{\mathbf j}$ for all $t\in T^\bullet_b$
in this case.
\begin{eqnarray}\label{BR}
g D_{\mathbf i} D_{\mathbf j} -
(g-\Lambda) D_{\mathbf j} D_{\mathbf i} & =& x_{\mathbf j} D_{\mathbf i}-
x_{\mathbf i} D_{\mathbf j}\ ,  \nonumber\\
g_s D_{\mathbf i} D_s-(g_s-\Lambda) D_s D_{\mathbf i} & = &
-x_{\mathbf i} D_s\,,
\quad \qquad \,\,\, \forall s\in S,
 \nonumber\\
g_s D_s D_{\mathbf j}-(g_s-\Lambda) D_{\mathbf j} D_s & = &
x_{\mathbf j} D_s \,, \qquad \qquad \,\, \forall s \in S,
 \nonumber\\
g_a^\circ :D_{\mathbf i} D_{t}: & = & -x_{\mathbf i} D_{t}\,,
\qquad \;\quad \,\, \forall t\in T_a^\circ ,\\
(g_a^\circ - \Lambda) :D_{\mathbf j} D_{t}: & = & -x_{\mathbf j} D_{t}\,,
\qquad \;\quad \,\, \forall t\in T_a^\circ ,  \nonumber\\
g_b^+ D_{\mathbf i} D_{t} & = & -x_{\mathbf i} D_{t}\,,
\qquad \;\quad \,\, \forall t\in T_b^\bullet ,  \nonumber\\
g_b^-  D_{t} D_{\mathbf j} & = & x_{\mathbf j} D_{t}\,,
\qquad \qquad \,\, \forall t\in T_b^\bullet , \nonumber
\end{eqnarray}
where  $g\neq 0$, $g_s\not=0$ for all $s$ and
$g_s\neq \Lambda$ for  $s$ such that either $s<{\mathbf i}$
or $s>{\mathbf j}$,
$g^\circ_a\not\in\lbrace  0, \Lambda\rbrace$ and $g_b^\pm\not = 0$.
\bigskip

\noindent
4.
\underline{Diffusion algebras of type $C(I=\{{\mathbf i}\},R_a)$:}
\bigskip

As in Lemma \ref{lem} the only element of $I$ is denoted here
as ${\mathbf i}$.
\begin{eqnarray}\label{C}
g_r D_{\mathbf i} D_r - (g_r- \Lambda_a) D_r D_{\mathbf i} & = & -
x_{\mathbf i} D_r \,,\quad   \forall r\in R_a ,
\end{eqnarray}
where $g_r\neq 0$ for $r<{\mathbf i}$ and
$g_r\neq \Lambda_a$ for $r>{\mathbf i}$.
\bigskip

\noindent
5.
\underline{Diffusion algebras of type $D(R)$:}

\begin{eqnarray}\label{lastOne}
D_r D_s - q_{sr} D_s D_r & = & 0 \,,\quad   \forall r<s\in R .
\label{lastOne}
\end{eqnarray}

\begin{thm}\label{Complete}
The list of algebras given above is exhaustive and
contains all possible Diffusion algebras with $N$ generators.
\end{thm}

\begin{proof}
According to the Diamond Lemma, an algebra of $N$ generators with relations
of type (\ref{algs}) is a Diffusion algebra if each of its triplets
$\{D_\alpha,D_\beta,D_\gamma\}$
generates a subalgebra coinciding with one of the cases listed in Section 2.
Lemma \ref{lem} provides information about possibile consistent
combinations of several such triplets and
we thus have to demonstrate that the families of algebras
(\ref{firstOne})--(\ref{lastOne}) exhaust the
list of Diffusion algebras which are allowed by this lemma.
\medskip

Let us start with the case $N_I\geq 3$.
According to  the first statement  of Lemma \ref{lem}
there are  two possible types of relations between generators from the set $I$.
This gives rise to  two families of Diffusion algebras --- $A_I$ and
$A_{II}$.
Statement 2  of Lemma \ref{lem} describes the relations between the
generators from the sets $I$ and $S$ in the case of the family $A_I$, and
the third statement of Lemma \ref{lem} excludes the presence of a nonempty set $S$
in the case of the family $A_{II}$.
The coefficients in the relations between the generators from the set $I$ and
the sets $T^\circ_a$ and $T^\bullet_b$
are subject to the conditions
given in the fourth statement of Lemma \ref{lem}, where
$\Lambda_{ij}=0$ for the $A_I$ family and $\Lambda_{ij}=(g_i-g_j)$
for the $A_{II}$ family.
These conditions fix the relations between the generators in the families of type
$A_I$ and $A_{II}$ to the expressions in (\ref{firstOne})
and (\ref{AII}).
Since the triplets of the form $\{I,S,T\}$ (occurring only for the $A_I$ family)
and $\{I,T_a,T_b\}$ for $a\not= b$
are of type $C^{(2)}$, no further conditions arise from these
relations. This exhausts all possibilities, and thus no further conditions
occur.
\medskip

Let $N_I=2$.
For each subset $S_a\subset S$
the connectivity property (\ref{extra2})
implies  relations of type $B^{(1)}$ for
all  triplets $\{I,I,S_a\}$, and thus
for any triplet $\{I,I,S\}$.
The corresponding relations are listed in the first three lines of (\ref{BR}).
The compatibility conditions within the triplets
$\{I,T_a,T_a\}$ and $\{I,I,T_a\}$ are given in the fourth statement of
Lemma \ref{lem}, where we now have
$\Lambda_{{\mathbf i}{\mathbf j}}\equiv\Lambda$.
These conditions fix the form of the last four lines in (\ref{BR}).
Since $C^{(2)}$ type relations
for $\{I,S,T\}$ and $\{I,T_a,T_b\}$ ($a\not= b$) triplets
do not imply further restrictions, no further
constraints arise.
\medskip

In the case $N_I=1$ a decomposition of the set $R$ into $S$ and
$T$ is not necessary, and we thus work with the whole set $R$.
Then the form of the relations (\ref{C}) is implied by the fifth statement of Lemma \ref{lem}
which describes the compatibility conditions for the $\{I,R_a,R_a\}$
triplets. The relations in the $C^{(2)}$ triplets $\{I,R_a,R_b\}$ for $a\neq b$
give no further constraints.
\medskip

In the case $N_I=0$ all the triplets are of type $D$,
which are compatible without any restrictions on the coefficients.
\end{proof}

Note that while mathematically possible, not all algebras in the families
are  relevant from the physicist's point of view. Due to the fact that the
structure constants of Diffusion algebras are interpreted as hopping rates,
that is probabilities, in the framework of stochastic processes on linear
lattices, only non-negative structure constants are relevant. This not only
implies restrictions on the structure constants themselves, but also on the
decompositions of the set $I_N$, because some configurations are not
compatible with non-negative structure constants. In particular, due to
Lemma \ref{lem}, statement 4 part 1, non-negative structure constants
throughout are possible only if the subsets $T^{\circ}_a$ fulfill one of the
following two requirements:
\begin{equation}\label{X}
1. \; \forall t\in T^\circ_a \mbox{~and~} \forall i\in I: t<i \quad \mbox{ or } \quad 2.
\; \forall t\in T^\circ_a \mbox{~and~}\forall i\in I: t>i\ .
\end{equation}

We conclude this section with some comments on the classification problem
for Diffusion algebras. To deal with the problem one should first establish
criteria of equivalence, and we  discuss two  natural ones here:
\medskip

\begin{itemize}
\item
One can consider linear transformations on the  set of generators
$\{D_\alpha\vert \alpha\in I_N\}$. However, there is the difficulty that not
all linear transformations respect the ansatz (\ref{algs}). There are two
special cases: rescaling transformations $D_\alpha\rightarrow \kappa_\alpha
D_\alpha$ and substitution transformations $D_\alpha \rightarrow
D_{\sigma(\alpha)}$, where $\sigma$ is an element of the symmetric group
${\cal S}_{N}$.

As has already been mentioned in the introduction,  rescalings may be used
to fix (depending on the context of the physical application)
some special values for the nonzero coefficients $x_\alpha$. In particular, this implies that the
values of the nonzero coefficients $x_\alpha$ are not relevant.

The substitution transformations clearly respect the form of the relations
(\ref{algs}), but may  contradict the requirement on the mutual ordering of
the generators, that is $g_{\alpha\beta}\neq 0$ for $\alpha<\beta$. In
particular, a permutation of the elements from different subsets $T_a$ and
between the subsets $T$ and $S$, or $T$ and $I$ is strictly forbidden. In
addition, one cannot permute two elements $r<s$ in the same subset $R_a$
unless $q_{sr}\not = 0$. On the other hand, permutations inside the subset
$I$ and (in most cases) between the subsets $I$ and $S$ are allowed unless
they contradict the reqirements described above. Thus, substitution
transformations establish certain equivalence classes inside each of the
families $A_I(I,S,T^\circ_a, T^\bullet_b)$,
$A_{II}(I,T^\circ_a, T^\bullet_b)$, $B(\{{\mathbf i},{\mathbf j}\},S,T^\circ_a, T^\bullet_b)$,
$C(\{{\mathbf i}\}, R_a)$ and $D(R)$. These equivalence classes  can be
calculated in concrete cases, but one hardly expects their complete
description in the case of general $N$.

Note that besides the rescalings and the substitutions which
do always exist there may occur other types of linear transformations which relate different
types of Diffusion algebras. For instance, in the case of $N=3$ the $C^{(1)}$ type
algebras in (\ref{C1}) with $\Lambda\neq 0$ can be reduced to  (a subclass of)
$D$ type algebras by the transformation $D_\alpha\rightarrow
D_\alpha-x_\alpha/\Lambda$. For general $N$, such transformations allow to
reduce  the number of nonzero parameters $\Lambda_a$ in the family of $C$
type Diffusion algebras in (\ref{C}) by 1.

\item
One can use the algebra antihomomorphism which inverts simultaneously multiplication
in the algebra, that is $D_\alpha D_\beta  \rightarrow D_\beta  D_\alpha$, and
the order of indices, that is $\alpha<\beta \rightarrow \alpha>\beta$.
This transformation amounts to a mirror reflection of the corresponding
stochastic processes. For example, in the list of Diffusion algebras with $N=3$
the families $B^{(3)}$ and $B^{(4)}$ are mirror symmetric. Further examples of
mirror symmetry for the case $N=4$ can be found in \cite{IPR} in
Appendix B.
\end{itemize}

\subsection{Description of the blending procedure}

The blending procedure is a constructive method to generate Diffusion
algebras. The corresponding Construction Theorem  states that any Diffusion
algebra can be obtained from a set of building blocks (equations
(4.1)--(4.7) in \cite{IPR}) via blending.
In the table below we describe the correspondence between the building blocks
from \cite{IPR} (left column) and
the specific subclasses of the families in Theorem \ref{Complete} (right column):
\begin{eqnarray}
A_I^{(1)} &:& A_I(I,S)\ ,\quad T=\emptyset\ ,\nonumber
\\
A_I^{(2)} &:&A_{I}(I,T^\circ)\ ,\;\; \mbox{and}\;\;  A_{II}(I,T^\bullet)\
,\quad S=\emptyset\ ,
\nonumber
\\
A_{II} &:& A_{II}(I,T^\circ)\ ,\;\; \mbox{and}\;\;  A_{II}(I,T^\bullet)\ ,
\label{Building} \\
B^{(1)}&:&
B(I=\{{\mathbf i},{\mathbf j}\},S)\ ,\quad T=\emptyset\ ,
\nonumber \\
B^{(2)} &:&
B(I=\{{\mathbf i},{\mathbf j}\},T^\circ)\ ,\;\; \mbox{and}\;\;
B(I=\{{\mathbf i},{\mathbf j}\},T^\bullet)\ ,\quad S=\emptyset\ ,
\nonumber \\
C &:&
C(I=\{{\mathbf i}\},R)\ ,  \nonumber \\
D &:& D(R)\ . \nonumber
\end{eqnarray}
Here it is understood that the sets $T^\circ$ and $T^\bullet$ whenever they
appear in the right
column of the table are the only
connective components in the decomposition (\ref{ST}) of the subset $T\in I_N$.
The connectivity condition may be also imposed on the subsets $S$ and $R$ in
the right column of (\ref{Building}).
We remark that the mathematial setting adopted in the present paper allow us to
extract elementary building blocks for the blending procedure. They are
shown in the right column of the table (\ref{Building}) and the blocks
listed in the left column and used in \cite{IPR} can be constructed by blending
of an arbitrary number of the corresponding blocks from the right column.
Note that  in the settings of Ref.\cite{IPR} extracting the elementary blocks
would only amount to imposing additional connectivity conditions (\ref{extra2}) on the
coefficients $g_{\alpha\beta}$ in the relations (4.1)--(4.7) there and so would
not suit the purposes of \cite{IPR}.

Furthermore we remark that
in contrast to \cite{IPR} we do not fix the order
between $D_{\mathbf i}$, $D_{\mathbf j}$ and $D_s$ to ${\mathbf i}<s<{\mathbf j}$
for all $s\in S$ in
$B(I=\{{\mathbf i},{\mathbf j}\},S)$, because the other orders
are needed when blending with $B(I=\{{\mathbf i},{\mathbf
j}\},T^\circ)$ and $B(I=\{{\mathbf i},{\mathbf
j}\},T^\bullet)$ in order to obtain all Diffusion algebras. This
is an inaccurracy in the formulation of the Construction Theorem in
\cite{IPR}. Despite that in the list of  $N=4$ Diffusion algebras given in
Appendix B
of \cite{IPR} the blending of such blocks
is treated correctly (see example 13 there).
\bigskip

Let $X_l(I, U_l)$, $l=1,\ldots,K$ denote $K$ building blocks in the list
(\ref{Building}) above, where $U_l$ refers to the set $R$,
$S$, $T^{\circ}$ or $T^\bullet$ corresponding to the building block,
and which are such that they have the same number of elements $N_I$ in $I$
with generators $D_i$, $i\in I$ satisfying in all blocks $X_l$ the same relations among
themselves.

Consider an ordered set $\cal I$ whose elements are labelled by the indices
from the sets $I$, $U_1,\dots ,U_K$
and such that for any $l=1,\dots ,K$
the order of the elements of $\cal I$ with their labels form $I$ and $U_l$
is the same as the order of the indices in the block $X_l(I,U_l)$. In this
situation we say that the order on $\cal I$ is {\em  compatible} with
the orders in the blocks $X_l$.

Let us denote as $X_{\cal I}(I, U_1,\dots , U_K)$ the algebra with generators labelled
by the elements in the set $\cal I$ and which satisfies the following conditions:
\begin{itemize}
\item
For any $l=1,\dots ,K$ the generators of $X_{\cal I}$ with indices from
the subsets $I, U_l\subset {\cal I}$ satisfy the same relations as their
corresponding generators in the blocks $X_l(I,U_l)$.
\item
For any $l_1\neq l_2\in\{1,\dots ,K\}$ and for all $a\in U_{l_1}\subset {\cal I}$
and $b\in U_{l_2}\subset {\cal I}$ the corresponding generators $D_a, D_b\in X_{\cal I}$
satisfy the relation
\begin{equation}
:D_a D_b: = 0\ .
\end{equation}
\end{itemize}
The procedure of constructing the algebras $X_{\cal I}$  from their building blocks
$X_l(I, U_l)$ is called {\em blending}.
Clearly the number of different algebras $X_{\cal I}$ which are associated
with the  set of building blocks $\{X_l(I, U_l)\}_{l=1,\dots ,K}$,
coincides with the number of different ordered sets $\cal I$ whose order is
compatible with the orders in all blocks $X_l$.

The following statement made in \cite{IPR} is a corollary to the Construction Theorem
\ref{Complete}:

\begin{thm}
Every Diffusion algebra can be obtained via a  blending of building blocks in
(\ref{Building}).
\end{thm}

\section{Conclusion}

We have presented a derivation of Diffusion algebras, which has led to five
different families of algebras: $A_{I}$, $A_{II}$ $B$, $C$ and $D$ and it
has been shown  that the approach is exhaustive. Since these families of
algebras correspond to the algebras obtained via the bending procedure in
\cite{IPR}, this also proves the Construction Theorem in this reference.

\section*{Acknowledgements}

We would like to thank A. Isaev and V. Rittenberg for useful discussions. R.T.
acknowledges furthermore useful discussions with A. Cox, V.K. Dobrev, P.P.
Martin and A. Sudbery.

\noindent
P.P. has partly  been supported  by the grant of  Heisenberg-Landau
Foundation  and by  RFBR grant \# 00-01-00299
and would like to acknowledge warm hospitality of
Max-Planck-Institut f\"{u}r Mathematik
where part of this work has been completed.
R.T. has partly been supported by a Marie Curie fellowship and is
grateful for the warm hospitality of the Department of Mathematics of the
University of York, where this work was started.

\end{document}